\magnification=\magstephalf
\vsize=7.5truein
\hsize=5.2truein

\loadeusm

\documentstyle {amsppt}

\def\slash{ /}
\def\a{ \alpha}
\def\b{ \beta}
\def\g{ \gamma}
\def\G{ \Gamma}
\def\d{ \delta}

\def\p{ \pi}

\def\si{ \sigma}
\def\Si{ \Sigma}

\def\bq{ \Bbb Q}
\def\br{ \Bbb R}
\def\bz{ \Bbb Z}

\def\8{_{ \infty}   }

\def\eR{ \eusm R}
\def\eS{ \eusm S}

\hfuzz 27pt

\topmatter \title  Topological Quantum Field Theory And Strong Shift
Equivalence  \endtitle
\rightheadtext{TQFT and SSE}    \author Patrick M.
Gilmer \endauthor
 \affil Louisiana State University \endaffil
 \address Department of Mathematics,  Baton Rouge, LA
70803 U.S.A  \endaddress
  \email gilmer\@ math.lsu.edu \endemail

\abstract Given a TQFT in dimension $d+1,$ and an infinite cyclic covering of a
closed $(d+1)$-dimensional manifold $M$, we define an invariant taking values in a strong
shift equivalence class of matrices. The notion of strong shift equivalence
originated in R. Williams'  work in symbolic dynamics. The Turaev-Viro module
associated to a TQFT and an infinite cyclic covering is then given by the Jordan
form of this matrix away from zero. This invariant is also defined if the
boundary of $M$ has a $S^1$ factor and the infinite cyclic cover of the boundary
is standard. We define a variant of a TQFT associated to a finite group $G$ which
has been studied by Quinn. In this way, we recover a link invariant due to
D. Silver and S. Williams. We also obtain a variation on the Silver-Williams
invariant, by using the TQFT associated to $G$ in its unmodified form. \endabstract
\thanks This research was partially supported by a grant from the Louisiana Education
Quality Support Fund \endthanks 
\keywords  Knot, Link, TQFT, Symbolic Dynamics, Shift
equivalence \endkeywords \subjclass 57R99, 57M99, 54H20 \endsubjclass
  \endtopmatter

\document 
\centerline  {This version:(1\slash 12\slash 98); 
First version: (7\slash 3\slash 97)}

\head Introduction\endhead

In this paper, we
will describe a relation between a knot invariant defined by D.Silver and
S.Williams \cite{SW1,SW2,SW3 SW4}, and Topological Quantum Field Theory (TQFT). As a
result, we will be able to give a new approach to the work of Silver and
Williams. The theorem of R. Williams which states that two shifts of finite type are conjugate if
and only if their adjacency matrices are strong shift equivalent (SSE) over the
nonnegative integers is replaced in our approach by the theorem which states that two
Seifert surfaces for the same link are stably equivalent. Both theorems have
elementary proofs. In this sense, the approach one might prefer depends on what
one already knows and feels comfortable with.

The connection that we are exploring between TQFT and symbolic dynamics is a two
way street. Borrowing the idea of strong shift equivalence from symbolic dynamics
allows us to define, given a TQFT, an invariant of an infinite cyclic covering 
that is in principle stronger than the associated Turaev-Viro module \cite{G1,G2}.   We have found a new invariant analogous to  the Silver-Williams invariant by using a TQFT
associated to a finite group as discussed by Quinn \cite{Q}.
On the other hand, Silver and Williams are finding further results in their
approach based on this connection. 

In the first section, we will define our SSE-class invariants. We do not actually
require all the axioms of a TQFT. We only need a functor from a (weak) cobordism category to the category of $R$-modules over some ring $r.$ In the second section, we
discuss Quinn's approach to a TQFT associated to a finite group. We also define
a variant of it. We discuss the first indications of the relation between these
TQFTs and the work of Silver and Williams.  In the third section, we combine the
results of sections one and two.  We  discuss the SSE-class invariants that we obtain
in this way, and the relation to the work of Silver and Williams.

We now give R. 
Williams' definition of SSE in its original form. Let $A, B$ be square matrices with nonnegative
entries. An elementary equivalence from $A$ to $B,$ is a pair of rectangular
matrices $(R,S)$ with nonnegative entries such that $A=RS$ and $B=SR.$ SSE is the
equivalence relation generated by elementary equivalences.\cite{W} \cite{LM}.
This is the definition of SSE that we will use in \S3.
In \S1, we will generalize this notion.

The initial data used by Silver and Williams is an ``augmented group system,''   while our  initial data is topological. Thus the results that Silver and Williams obtain apply 
 directly to questions in combinatorial group theory.  The connection we found between their work and TQFT has led them to study the shift space we describe at the end of \S3. They have also found another shift of finite type which allows them to find a linear recursion formula for the number conjugacy classes of representations of the fundamental group of the branched cyclic covers of a knot.  A future direction suggested by \cite{Q, \S 5} is the introduction of cohomology classes for the finite group $G$ into the initial data.  
Also needing further exploration is the use of the more general notion of SSE introduced in \S 1 to the study of, for instance, the TQFT's  discused in \cite{BHMV}.

For convenience we work with smooth oriented manifolds.

We thank Dan Silver and Susan Williams for useful discussions. 

\head \S 1 SSE-class invariants\endhead

We will consider pairs $(M,\chi),$ where $M$ is a compact connected 
$(d+1)$-manifold  and $\chi$ is a primitive element  in $H^1(M).$ 
Alternatively, $\chi$ is an epimorphism from $H_1(M)$ to $\bz.$  We can think of
$\chi$ as a choice of a connected infinite cyclic covering space (together with a choice of generator for the group of covering translations). If $M$ has
boundary we assume, that the boundary has the form $P\times S^1,$ for some
$(d-1)$-manifold $P.$ Moreover we insist that $\chi$ restricted to the boundary
is given by the projection on the second factor under the canonical isomorphism
$H_1(\partial M)= H_1(P\times S^1)\approx H_1(P)\oplus \bz.$ Thus the  induced
cover of $P\times S^1$ is $P \times \br,$ with the standard action. We may
describe the situation that $M$ is without boundary by saying that $P$ is empty.

A Seifert manifold for $(M,\chi)$ is a connected  $d$-manifold, $\Si,$
Poincare-Lefshetz dual to $\chi \in H^1(M),$  such that $\Si \cap \partial M
=\partial \Si= P \times {z},$ for some $z\in S^1.$

 We consider the following type of elementary modification of a Seifert manifold
$\Si.$  Suppose we have an embedding $f$ of $D^k \times D^{d-k+1}$ in the interior
of $M$ and $f(D^k \times D^{d-k+1})\cap \Si= f(S^{k-1} \times D^{d-k+1}).$ We
assume here that $k \le \frac {d+1} 2.$  Then we may let $\Si'$ denote
$(\Si-f(S^{k-1} \times D^{d-k+1}))\cup f (D^k \times S^{d-k})$ after smoothing.
$\Si'$ is also a Seifert manifold. We will say $\Si'$ has been obtained by
elementary expansion of $\Si.$ The equivalence relation on Seifert manifolds
generated by elementary expansions and isotopies is called stable equivalence.

Let $L$ be an oriented link in $S^3.$ Then we may take $M$ to be the exterior of
$K$ together with $\chi$ given by the linking number with $L$ in $S^3.$ Then a
Seifert manifold for $(M,\chi)$ is a Seifert surface for $L$ in the classical
sense.

We have the following generalization of the stable equivalence of Seifert
surfaces for links.  It may be proved by adapting the proof given in
\cite{G-L, page 66}. The case where $M$ is a classical knot exterior  has another
elementary proof \cite{BF}.

\proclaim{Lemma 1.1} Let $(M,\chi)$ be as above, then any two Seifert manifolds for
$(M,\chi)$ are stably equivalent.\endproclaim

Given a Seifert manifold  $\Si$ for $M,$ we may consider, $E_\Si,$ $M$ slit along $\Si.$
We may view this as a cobordism rel boundary from $-\Si$ to $\Si.$  Stacking a
bi-infinite sequence of copies of $E_\Si,$ we obtain $M_\chi$ the infinite cyclic
cover of $M$ classified by $\chi.$ We note $E_\Si$ is a fundamental domain for the
$\bz$ action on $M_\chi.$ 

By a weak cobordism category, we mean the usual notion of cobordism category except the operations of disjoint union need not be defined on either objects or morphisms.  See \cite{BHMV} or
\cite{Q} for the notion of cobordism category.
Let $\Cal C(P)$ denote the weak cobordism category whose objects consist 
of  connected 
$d$-manifolds with boundary $P$ and whose morphisms from object $\Si$ to object $\Si'$
are (equivalence classes of) connected cobordisms rel boundary from $-\Si$ to $\Si'.$  As usual, two connected cobordisms rel boundary from $-\Si$ to $\Si'$ are 
equivalent if there is a diffeomorphism between them respecting the identification of parts of their boundary to $-\Si \coprod \Si'$

Let $r$ be a commutative ring. Now suppose that $(Z,V)$ is a
functor from $\Cal C(P)$ to the category of finitely generated $r$-Modules, and
$r$-Module homomorphisms. We follow the tradition in TQFT of letting $Z$ denote
the application of the functor to morphisms and letting $V$ denote the application of
the functor to objects. Let $Z(\Cal C(P))$ denote the algebroid  consisting of all
the  $r$-Module homomorphisms obtained by applying $Z$ to a morphism of $\Cal C(P).$

Given $A, B \in Z(\Cal C(P)),$ an elementary equivalence from $A$ to $B$ is a
pair of elements $R,S \in Z(\Cal C(P))$ such that $A=RS$ and $B=SR.$ We define strong shift equivalence (SSE) to be the equivalence
relation on $Z(\Cal C(P))$ generated by elementary equivalences.

Given  $A_1:V(\Si) \rightarrow V(\Si),$ and $A_2: V(\Si) \rightarrow V(\Si),$ we will say $A_1$
and $A_2$ are similar if there exists a diffeomorphism $f:\Si \rightarrow \Si$
such that $Z(f)A_1 = A_2 Z(f).$ Here $Z(f)$ is the map given by $Z(C_f)$ where
$C_f$ is  the mapping cylinder of $f.$  A similarity is an elementary equivalence: let $R= A_1 Z(f^{-1}),$ and $S= Z(f),$ then $RS=A_1$ and $SR=A_2.$
Thus changing our identification of $\Si$ with a model surface does not change the SSE class of $Z(E_{\Si}).$

If $\Si'$ is an elementary expansion of $\Si,$ let $\eR$ be formed by attaching a
k-handle $D^k \times D^{d-k+1}$ to $\Si \times I$ using the embedding $f$
restricted to $S^{k-1} \times D^{d-k+1}$ as an attaching map. The $\eR$ embeds
in $M$ so it is a relative cobordism between $\Si$ and a pushed off copy of
$\Si'.$ Let $\eS$ denote the closure of $E_\Si-\eR.$ Then $\eS \circ \eR = E_\Si$
and $\eR \circ \eS = E_{\Si'}.$ Thus $(Z(\eR),Z(\eS))$ is an elementary
equivalence between $Z(E_\Si)$ and $Z(E_{\Si'}).$ 

Moreover an isotopy between 
$\Si$ and $\Si'$  induces compatible diffeomorphisms between $\Si$ and $\Si'$ and $E_\Si$ and $E_{\Si'}$ which yield an elementary
equivalence between $Z(E_\Si)$ and $Z(E_{\Si'}).$

\proclaim{Theorem 1.2} For any two Seifert manifolds $\Si$ and $\Si',$
 $Z(E_\Si)$ and $Z(E_{\Si'})$ are SSE in  $Z(\Cal C(P)).$ \endproclaim

We often may need to consider such extra structure on the manifolds used as
objects and morphisms of $\Cal C (S)$ as a $p_1$-structure or a submanifold (of possibly fixed
codimension). The above arguments easily adapt. The Seifert manifolds should be
taken in general position with respect to the submanifold.

In practise, we may not know much about $Z(\Cal C(P)).$ Then we may replace $Z(\Cal
C(P))$ by a larger algebroid, for instance all $r$-module homomorphisms.  

If $r$ is an integral domain, we may tensor with the field of fractions, and then may take the similarity class of the invertible
part of $Z(E_\Si).$  This is essentially the
Turaev-Viro module $Z(M,\chi)$ discussed in \cite{G1} (in the case $P=
\emptyset$).  This procedure is discussed in \cite{LM,7.4} in the case when $r=\Bbb Z,$  and the modules are free.  
There is a weaker notion of equivalence than strong shift
equivalence called shift equivalence.  Shift equivalence is easier to analyze
algebraically \cite{LM,7.5}.   $Z(M,\chi)$ is determined by the shift equivalence
class of $Z(E_\Si).$ An example over $r=\Bbb Z$  given in \cite{LM,7.3.4} shows that
shift equivalence is a finer invariant than the similarity class of the invertible part over $\bq.$

For the functors that we  consider in the next section, $V(\Si)$ is always a free
$\bq$-vector space or $\bz$-module with a fixed unordered basis. Diffeomorphisms induce maps on
$V(\Si)$ which preserve this unordered basis. Moreover matrices which represent
elements of $Z(\Cal C(P))$ with respect to the given bases always have nonnegative
integral entries.   For this, it is important that the objects of $\Cal C(P)$ be connected.  

The above theorem then implies that matrices which represent $Z(E_\Si)$ and
$Z(E_{\Si'})$ with respect to these bases are SSE in the original sense of
R. Williams.

\head \S 2 Two TQFTs associated to a finite group $G$\endhead

Quinn discusses a TQFT associated to a finite group $(Z_G,V_G).$  Special cases
were earlier studied  by Kontsevich, Dijkgraaf-Witten, Segal and Freed-Quinn \cite{FQ}. We will generally follow Quinn's development \cite{Q}. However we
depart from Quinn's use of a single letter $Z$ for the application of the functor on
objects and morphisms.

Quinn discusses how a cobordism category is constructed from a ``domain category''
which consists of a pair of categories: ``spacetimes'' and  ``boundaries''.
$(Z_G,V_G)$ is defined on a cobordism category $\Cal C$ whose objects
(``boundaries'') are  finite CW complexes.  A morphism  from $Y_1$ to $Y_2$ is
an equivalence class  of finite CW pairs $(X,Y)$ together with a
homeomorphism $f:Y \rightarrow Y_1 \coprod Y_2$ of boundaries.  Two such triples
$(X,Y,f),$ and  $(X',Y',f')$ are equivalent if there is a homotopy equivalence
$G$ from $(X,Y)$ to $(X',Y')$ which makes the following diagram \cite{Q,p. 369}  
commute.

$$\CD
 Y @> G_{|Y} >> Y'\\
 @V{f}VV         @VV{f'}V  \\
 Y_1 \coprod Y_2 @> \text{identity} >> Y_1 \coprod Y_2 \endCD $$ 

The vector space associated to an object $Y$ is the rational vector space with
basis $[Y,BG].$  If $Y$ is connected, this can be identified with
$\hom(\p_1(Y),G)\slash G,$ the representations of the
fundamental group of $Y$ into the group $G,$ modulo the action given by
conjugation by $G.$ Suppose $(X,Y)$ is a finite CW pair, and $Y$ is provided with
a homeomorphism to
 $Y_1 \coprod Y_2,$ where $Y_1$ is to be thought of as the incoming boundary, and
$Y_2$ is to be thought of as the outgoing boundary. Then $Z(X,Y):V(Y_1)
\rightarrow V(Y_2)$ is defined. The definition of this map in down to earth terms
may be given in various cases in terms of the number of components of $Y_1$
\cite{Q,4.13}.  We only mention the two formulas that we will use. If $Y$ is
empty, so $X$ is closed, then $Z(X)$ is an automorphism of $V(\emptyset)$ which
is one dimensional. Thus $Z(X)$ is essentially a rational number. If $X$ is connected and $Y$ is empty, we have:

$$Z_G(X)= \sum_{[\b] \in \hom(\p_1(X),G)\slash G} \frac {1}{\# C_\b}= \sum_{[\b]
\in \hom(\p_1(X),G)\slash G} \frac {\#[\b]}{\# G}= \frac {\#\hom(\p_1(X),G)}
{\#G}\leqno {2.1}$$

\noindent Here we use $\#$ to denote the number of elements in a set, and $C_\b$ denotes
the centralizer in $G$ of the image of $\b.$  This is taken from \cite{FQ,5.14}
where we correct a typo.

Suppose $X,$ $Y_1,$ and $Y_2$ are all connected. Let us pick a path $\si$ in $X$
from $y_1 \in Y_1$ to $y_2 \in Y_2.$ $\si$ induces a homomorphism
$\si_*:\pi_1(X,y_2) \rightarrow \pi_1(X,y_1),$ which sends a loop $\g$ at $y_2$
to the loop $\si\g \si^{-1}$ at $y_1.$ We may identify $V(Y_i)$ with
$Q[\hom(\pi_1(Y,y_i),G)\slash G].$ If $\b:\pi_1(X,y_1) \rightarrow G,$ let $\b_1:
\pi_1(Y_1,y_1) \rightarrow G$ denote the map induced by inclusion from 
$\pi_1(Y_1,y_1) \rightarrow \pi_1(X,y_1)$ composed with $\b.$ Let $\b_2:
\pi_1(Y_2,y_2) \rightarrow G$ denote the map induced by inclusion from 
$\pi_1(Y_2,y_2) \rightarrow \pi_1(X,y_2)$ composed with $\si_*$ followed by $\b.$
For simplicity we are pretending that the boundary identification $f$ is the
identity in these formulas. This should cause no confusion. If $\a \in
\hom(\pi_1(Y_1,y_1),G)$ 
$$Z_G(X,Y)([\a])= \sum_{\b \in \hom(\pi_1(X,y_1) , G) \ni
\b_1=\a} [\b_2]\leqno {2.2}$$

Given a morphism $(X,Y,f)$ from $Y_1$ to itself, we may form the mapping torus
$T((X,Y,f))$ by identifying the two components of $Y$ in $X$ using $f.$ By \cite{Q,7.5},
we have that $Z_G(T((X,Y,f)))$ is given by the trace of $Z_G(X,Y,f).$

Let $(M,\chi)$ be as in \S 1, and $\Si$ a Seifert manifold for $(M,\chi).$  Let $M_d$ denote the associated $d$-fold cyclic
cover of $M.$ Let $E_{\Si}$ now denote the morphism from $\Si$ to itself given by
$E_{\Si}.$ Note that $M_d= T( E_{\Si}^d).$ The Cayley-Hamilton Theorem implies that the powers of a matrix satisfy a linear recursion relation. Taking the trace of this 
relation,  we obtain \cite{G1,1.8}

\proclaim {Proposition 2.3} $\# \hom(\p_1(M_d),G)$ satisfies a linear recursion
formula in $d$ whose coefficients are given by the characteristic polynomial of
$Z_G(E_{\Si}).$ \endproclaim

We wish now to define $(\hat Z_G,\hat V_G),$ a variation of the above TQFT. We
must first describe a new domain category. The spacetimes will be triples
$(X,Y,W)$ where $(X,Y)$ is a finite CW pair, and $W$ is a subcomplex of
$X$ such that the intersection of $W$ with each path component of $Y$ and each path component of $W$ is
nonempty.  An object is a ``boundary,'' i.e.\  a CW pair $(Y,U)$ where each path
component of $Y$ has a nonempty intersection with $U.$  $\hat V(Y,U)$ is the
free $\bz$-module  with basis $[(Y,U),(BG,*)].$  This basis can be identified
with $\hom(\p_1(Y\slash U),G),$ in the case that $Y$ is nonempty.   A
morphism from $(Y_1,U_1)$ to $(Y_2,U_2)$ is now an equivalence class of a 4-tuple
$(X,Y,W,f)$ where   $f:Y \rightarrow Y_1 \coprod Y_2 $ is a homeomorphism and 
$U_i$ denotes $f(W \cap Y)\cap Y_i.$  The equivalence relation is the same as
above for $(X,Y),$ except that no longer we require that the map   $G:X
\rightarrow X'$ be a homotopy equivalence but we do require that $G$ map $W$ to $W',$ and  that $G$ induces a homotopy equivalence
from $X\slash W$  to $X'\slash W'.$  We let $[W]$ denote the point in $X\slash W$
which is the image of $W$ under the quotient map. Similarly for $[U_i]$ etc.
One may check that Quinn's axioms for  a domain category are satisfied.

We now define  $\hat Z_G(X,Y,W)$ in elementary terms. Note that we do not need
 to 
take into consideration the number of components of $Y_1.$
If $\d:\pi_1(X\slash W,[W]) \rightarrow G,$ for $i$ equal to one or two, let $\d_i:
\pi_1(Y_i,U_i) \rightarrow G$ denote the map induced by inclusion from 
$\pi_1(Y_i\slash U_i,[U_i]) \rightarrow \pi_1(X\slash W,[W])$ composed with $\d.$
We are still pretending that $f$ is the identity in these formulas. We define

$$\hat Z_G(X,Y,W)(\g) = \sum_{\d \in \hom(\pi_1(X\slash W,[W]) , G) \ni \d_1=\g} \d_2 \leqno {2.4}$$

The functoriality of this map with respect to composition of morphisms now
follows directly from  the Van-Kampen theorem as discussed for instance in
Massey \cite{M}. If $Y$ is empty, we have:

$$\hat Z_G(X,W)= \# \hom(\pi_1(X\slash W,[W]) , G).\leqno {2.5}$$

Let $(M,\chi)$ be as in \S 1. Suppose $M$ contains a submanifold $N$ in its
interior. Suppose that $\chi$ restricted to $N$ is nontrivial.  If $\Si$ is a Seifert manifold transverse to $N,$ then $N \cap \Si$ is nonempty. Let $N_d$ denote
the associated $d$-fold cyclic cover of $N.$ Note we may view $N_d$ as a
submanifold of $M_d.$  Let $(E_{\Si},N \cap E_{\Si})$ denote the morphism from
$(\Si,\Si\cap N)$ to itself given by $(E_{\Si},N \cap E_{\Si}).$ We have the following analog of (2.3).

\proclaim {Proposition 2.6} $\# \hom(\p_1(M_d\slash N_d),G)$ satisfies a linear
recursion formula in $d$ whose coefficients are given by the characteristic
polynomial of $\hat Z_G(E_{\Si},N \cap E_{\Si}).$\endproclaim

Suppose $L$ is a oriented codimension-2 link in $S^n$ with $\mu$ components, $M$ is its  exterior, $\chi$ the map which sends  each oriented meridian to one, and $N$
is the disjoint union of one meridian for each component of $L.$  Let  $B_d$ denote the $d$-fold cyclic branched cover of $S^n$ along $L.$ 
Let $Q_d$ denote $M_d$ with  $\mu$ 2-disks attached along the components of $N_d.$ Then $M_d \slash N_d$ is homotopy equivalent to the space obtained from $Q_d$ by identifying the  $\mu$ center points of these added 2-disks. Up homotopy,  this is equivalent to joining these points with $\mu-1$ arcs.
So $\pi_1(M_d \slash N_d) = F_{\mu-1} *\pi_1(Q_d).$ If $n=2,$ $B_d=Q_d.$ If
$n>2,$ $B_d$ may be obtained from $Q_d$ by thickening the $\mu$
2-disks and then adding $\mu$ n-disks along their whole boundaries. Thus
$\pi_1(B_d)= \pi(Q_d).$ Thus $\pi_1(M_d \slash N_d)=\p_1(B_d) * F_{\mu-1}.$
Since  $\# \hom(\p_1(B_d) * F_{\mu-1},G) = (\#G)^{\mu-1}\left(\# \hom(\p_1(B_d) ,G) \right),$ $\# \hom(\p_1(B_d) ,G )$ also satisfies a linear recursion formula whose coefficients are given by the characteristic
polynomial of $\hat Z_G(E_{\Si},N \cap E_{\Si}).$ Thus the above
proposition generalizes \cite{SW2,Corollary (4.2)}. This observation, in the case of knots, was the second
indication that there was a connection between \cite{G1}, and \cite{SW2}.

The
first indication of a connection was independent derivations of the following result by W.
Stevens \cite{S1,S2}, and Silver-Williams \cite{SW2}. For any finite abelian
group $G,$ the homology of the $d$-fold branched cyclic cover of $S^3$ along a fixed
knot with coefficients in $G$ is periodic in $d.$ I had conjectured this result
based on computer experiments based on
Proposition (2.3)  with $G=Z_{p^r}.$  I gave this as a thesis problem to my student Stevens. He gave an nice proof using classical techniques. Silver and Williams independently discovered this result using symbolic dynamics.

We conclude this section with a remark which we have not seen in the literature.
The TQFT's discussed in this section
are not generally ``cobordism generated'' in the sense of \cite{BHMV} or
what is the same thing ``nondegenerate'' in the sense of \cite{T,III.3.1}. To see this consider $G=Z_p,$ with $p$ a prime, then $Z_G(Y)$ is $\bq[H_1(Y,Z_p)].$
A vacuum state must have constant coefficients along the rays of $H_1(Y,Z_p).$
By ray, we mean a line through the origin with the origin deleted. 

\head \S 3 Symbolic dynamics \endhead

Suppose $M$ is a $(d+1)$ dimensional manifold with $N$, a submanifold.
We are especially thinking of the case $M$ is either a link exterior and $N$
is a meridian for a chosen component of the link, or $M$ a
homology $S^d \times S^1 $ obtained by performing surgery along a knot, and $N$ is empty. There are, of course, many other  possibilities. Let $\chi$ be a primitive cohomology class
as in \S1 which remains nontrivial when pulled back to $N.$  Let $\Si$ be
a Seifert manifold transverse to $N,$  and  $E_{\Si}$ be $M$ slit along $\Si.$ 

We may apply $(\hat Z_G,\hat V_G)$ to the endomorphism $(E_{\Si}, N \cap E_{\Si})$ 
of $({\Si}, N \cap {\Si}).$ $\hat V_G({\Si}, N \cap {\Si})$ is a free $\bz$-module which comes equipped with a choice of unordered basis. With respect to this choice of basis and any ordering, a matrix for
$\hat Z_G(E_{\Si}, N \cap E_{\Si})$ has entries which are non-negative integers.

Thus by \S 1, the
SSE equivalence class of a matrix for   $\hat Z_G(E_{\Si}, N \cap E_{\Si})$ with respect to this basis
is an invariant of the triple $(M,N,\chi),$ which we denote by
$I_G(M,N,\chi).$

Given the endomorphism
$(E_{\Si},E_{\Si}\cap N),$  we may
associate a directed graph $\hat \G$ as follows. $\hat \G$ has  one vertex for each basis
element for $\hat V_G(\Si,\Si \cap N),$ i.e.\ each element of 
$\hom(\pi_1(\Si\slash ( {\Si}\cap N)),G).$  We draw $a(v,v')$ edges from 
vertex $v$ to vertex $v'$ if $a(v,v')$ is the $v'$ coefficient of
$Z_G(E_{\Si}(v)).$   Thus the matrix for $\hat Z_G(E_{\Si})$ is the adjacency matrix
for $\hat \G.$  

According to the classification theorem of R. Williams \cite{LM,7.2.7},   topological conjugacy of the
shift space $\Cal X_{\Cal G}$ of bi-infinite paths in a directed graph $\Cal G$ corresponds
exactly to  SSE equivalence of the adjacency matrix for $\Cal G.$ 
Consider the shift space $\Cal X_{\hat \G}$ of finite type given by all
bi-infinite paths in $\hat \Gamma.$  By the classification  theorem, the topological equivalence
class of $\Cal X_{\hat \G}$ is also invariant of  the triple $(M,N,\chi).$  This
invariant contains the same information as $I_G(M,N,\chi).$

Note that the edges of $\hat \G$ can be put in 1-1 correspondence with 
$\hom(\pi_1(E_{\Si}\slash (E_{\Si}\cap N),G).$  The edge corresponding to
$\rho \in \hom(\pi_1(E_{\Si}\slash (E_{\Si}\cap N),G)$ goes from the vertex
corresponding to the ``restriction'' of  $\rho$ to the incoming copy of $\Si$ to the
vertex which corresponds to $\rho$ restricted to the outgoing  copy of $\Si.$ 
See (2.4). Thus  a bi-infinite path in  $\Cal X_{\hat \G}$ can be used to define a representation from $\pi_1(M_\chi\slash N_\chi,[N_\chi])$ to $G$ and visa-versa.
Here $M_\chi$ denotes the infinite cyclic covering of $M$ given by $\chi.$
$N_\chi$ is the associated infinite cyclic covering of $N.$ Both of these spaces and their quotients have $\bz$ actions on them. Thus $\hom(\pi_1(M_\chi\slash N_\chi,[N_\chi]),G)$ is the shift space which corresponds to $\Cal X_{\hat \G}.$

If $M$ is a  knot or link exterior, and $N$ is a meridian to the chosen component, this can be seen to be the shift space studied by Silver and Williams \cite{SW1,SW2}.

Now we may also apply the functor $(Z_G,V_G)$ to the endomorphism
$E_\Si$ of $\Si$ where $N$ is now taken to be empty. $V_G({\Si})$ is a rational vector space which comes equipped with a choice of unordered basis. Recall that we specified that Seifert manifolds must be connected in \S 1. It follows that  
with respect to the given choice of basis and any ordering, a matrix for
$Z_G(E_{\Si})$ has entries which are non-negative integers.

By \S 1, the
SSE equivalence class of a matrix for   $Z_G(E_{\Si})$ with respect to this basis
is an invariant of the pair $(M,\chi),$ which we denote by
$I_G(M,\chi).$
 
As above we may also describe this new invariant as the topological conjugacy class of the shift space  given by bi-infinite paths in a directed graph $\G$ whose adjacency matrix is the matrix for $Z_G(E_{\Si}).$

$\G$ has  one vertex for each element of $\hom(\pi_1(\Si,G)\slash G.$
In order to describe the edges of $\G$ it is useful to consider the directed graph $\hat \G$ described above, but now assuming the $N$ is a circle which meets $\Si$ in a single point. Then $G$ acts by conjugation of the vertices of $\hat \G,$ 
and the vertices of $\G$ are the orbits of this action.
If $v$ stands for a vertex of $\hat \G$, let [v] denote the vertex of $\G$ given by the orbit of $v.$  
Comparing (2.2) and (2.4), we see that the number of edges joining $[v]$ to $[w]$ in $\G$ is the number of edges in $\hat \G$ joining $v$ to some vertex in the orbit of $w.$

  \enddocument \vfill\eject \end